\documentstyle[12pt,amssymb,amscd]{amsart}  

\newtheorem{thm}{Theorem}[subsection]

\newtheorem{cor}[thm]{Corollary}
\newtheorem{prop}[thm]{Proposition}
\newtheorem{lemma}[thm]{Lemma}

\newtheorem{conj}[thm]{Conjecture}

\theoremstyle{remark}
\newtheorem{remark}[thm]{Remark}

\theoremstyle{definition}

\numberwithin{equation}{section}

\newcommand{\bbA}{{\Bbb A}}
\newcommand{\bbC}{{\Bbb C}}

\newcommand{\bbL}{{\Bbb L}}
\newcommand{\bbR}{{\Bbb R}}
\newcommand{\bbW}{{\Bbb W}}
\newcommand{\bbZ}{{\Bbb Z}}
\newcommand{\cF}{{\cal F}}

\newcommand{\cO}{{\cal O}}

\newcommand{\cM}{{\cal M}}
\newcommand{\cN}{{\cal N}}
\newcommand{\cD}{{\cal D}}

\newcommand{\cA}{{\cal A}}

\newcommand{\cI}{{\cal I}}

\newcommand{\cC}{{\cal C}}
\newcommand{\cE}{{\cal E}}

\newcommand{\cU}{{\cal U}}
\newcommand{\cR}{{\cal R}}
\newcommand{\cS}{{\cal S}}

\newcommand{\chr}{\operatorname{char}}
\newcommand{\supp}{\operatorname{Supp}}
\newcommand{\SS}{\operatorname{SS}}

\newcommand{\Eu}{\operatorname{Eu}}
\newcommand{\eu}{\operatorname{eu}}

\newcommand{\R}{\operatorname{\bold R}}
\newcommand{\Hom}{\operatorname{Hom}}
\newcommand{\isomo}{\overset{\sim}{=}}

\newcommand{\id}{\operatorname{id}}

\newcommand{\shHom}{\underline{\operatorname{Hom}}}

\newcommand{\ad}{\operatorname{ad}}
\thanks
{The third author was supported
in part by NSF grant DMS-9504522.}

\begin{document}

\title{Riemann-Roch theorems via deformation quantization I}

\author{P.Bressler}
\address{Department of Mathematics, The Pennsylvania State University,
University Park, PA 16802, USA}
\email{bressler@@math.psu.edu}
\author{R.Nest}
\address{Mathematical Institute, Universitatsparken, 5,
2100 Copenhagen, Denmark}
\email{rnest@@math.ku.dk}
\author{B.Tsygan}
\address{Department of Mathematics, The Pennsylvania State University,
University Park, PA 16802, USA}
\email{tsygan@@math.psu.edu}

\maketitle

\begin{abstract}
We deduce the Riemann-Roch type formula expressing the microlocal Euler
class of a perfect complex of $\cD$-modules in terms of the Chern character
of the  associated symbol complex and the Todd class of the manifold from
the Riemann-Roch type theorem for periodic cyclic cocycles of a symplectic
deformation quantization. The proof of the latter is contained in the sequel
to this paper.
\end{abstract}

\section{Introduction}
We recall the results of \cite{SS} restricting ourselves to the ``absolute''
case for the sake of simplicity.

\subsection{The index theory for elliptic pairs}\label{section:ell}
Let $X$ be a complex manifold of dimension $\dim_{\bbC}X = d$.
Say that a $\cD_X$-module $\cM$ is {\em good} is for every relatively
compact open subset $U$ of $X$ the restriction $\cM\vert_U$ admits a finite
filtration $G_\bullet\cM\vert_U$ by $\cD_U$-submodules, such that
$Gr^G_\bullet\cM\vert_U$ is $\cD_U$-coherent with good filtration.

An {\em elliptic pair} $(\cM^\bullet,\cF^\bullet)$ on $X$ consists
of
\begin{itemize}
\item a complex $\cM^\bullet$ of $\cD_X$-modules with bounded good cohomology,
\item a complex $\cF^\bullet$ of $\bbC$-vector spaces with bounded,
$\bbR$-constructible cohomology,
\end{itemize}
which satisfy
\[
\chr(\cM^\bullet)\cap\SS(\cF^\bullet)\subseteq T^*_XX\ .
\]

In the case when $\supp\cM^\bullet\cap\supp\cF^\bullet$
(where the support of a complex of sheaves is understood to be the
union of the supports of the cohomologies) is compact, P.~Schapira and
J.-P.~Schneiders proved that
\[
\dim H^\bullet(X;\cF^\bullet\otimes\cM^\bullet\otimes^{\bbL}_{\cD_X}\cO_X)
<\infty \ .
\]
Thus, the Euler characteristic
\[
\chi(X;(\cM^\bullet,\cF^\bullet))\overset{def}{=}\sum_{i} (-1)^i
\dim H^i(X;\cF^\bullet\otimes\cM^\bullet\otimes^{\bbL}_{\cD_X}\cO_X)
\]
is defined. The index theorem for elliptic pairs (\cite{SS}, Theorem 5.1)
says that 
\[
\chi(X;(\cM^\bullet,\cF^\bullet))=
\int_{T^*X}\mu\eu(\cM^\bullet)\smile\mu\eu(F^\bullet)
\]
where $\mu\eu(\cM^\bullet)\in H^{2d}_{\chr(\cM^\bullet)}(T^*X;\bbC)$ and 
$\mu\eu(F^\bullet)\in H^{2d}_{\SS(F^\bullet)}(T^*X;\bbC)$ are defined in
\cite{SS}.

The class $\mu\eu(F^\bullet)$ is the characteristic cycle of the constructible
complex $F^\bullet$ as defined by Kashiwara (see \cite{KS} for more details).
For example, if $Y\subset X$ is a closed real analytic submanifold, one has
$\mu\eu(\bbC_Y) = [T^*_YX]$.

With regard to $\mu\eu(\cM^\bullet)$ P.~Schapira
and J.-P.~Schneiders conjectured that it is related to a certain 
characteristic class of the symbol $\sigma(\cM^\bullet)$ of $\cM^\bullet$
(\cite{SS}, Conjecture 8.5, see Conjecture \ref{conj:ch}).

\subsection{The Riemann-Roch type formula}
Let $\pi : T^*X\to X$ denote the canonical projection.

In what follows we will assume that the complex $\cM^\bullet$ admits a
filtration $F_\bullet\cM^\bullet$ by $\cO_X$-submodules which is compatible
with the action of $\cD_X$ and the filtration $F_\bullet\cD_X$ by order and
such that the {\em symbol complex}
of $\cM$ defined by
\[
\sigma(\cM^\bullet) = \pi^{-1}Gr^F_\bullet\cM^\bullet\otimes_{\pi^{-1}
Gr^F_\bullet\cD_X}\cO_{T^*X}
\]
has bounded $\cO_{T^*X}$-coherent cohomology. Note that, by definition,
$\chr(\cM^\bullet)\overset{def}{=}\supp\sigma(\cM^\bullet)$.

For $\Lambda$ a closed subvariety of $T^*X$ let
$K^0_\Lambda(T^*X)$ denote the Grothendieck group of perfect complexes of
$\cO_{T^*X}$-modules supported on $\Lambda$ (i.e. acyclic on the complement
of $\Lambda$ in $T^*X$). Let $K^0(T^*X)\overset{def}{=}K^0_{T^*X}(T^*X)$.

For $\Lambda$ containing $\chr(\cM)$ let 
$\sigma_\Lambda(\cM^\bullet)$ denote the class of $\sigma(\cM^\bullet)$
in $K^0_\Lambda(T^*X)$.

\begin{remark}
Both the characteristic variety and the class of the symbol in the
Grothendieck group are independent of the choice of the good filtration.
Local existence of good filtration in coherent $\cD_X$-modules is sufficient
to define $\sigma_\Lambda(\cM^\bullet)$ for $\cM^\bullet$ with bounded
good cohomology.
\end{remark}

The Chern character $ch : K^0(T^*X)\to\bigoplus_i H^{2i}(T^*X;\bbC)$ admits
a natural extension to the Chern character with supports
(Theorem \ref{thm-intro:char-cl})
\begin{equation}\label{map:ch-with-supp}
ch_\Lambda : K^0_\Lambda(T^*X)\to\bigoplus_i H^{2i}_\Lambda(T^*X;\bbC)
\end{equation}
which is functorial with respect to change of support.

In \cite{SS}, P.~Schapira and J.-P.~Schneiders make the following conjecture.
\begin{conj}\label{conj:ch}
Suppose that $\cM^\bullet$ is a complex of $\cD_X$-modules with bounded good
cohomology and $\Lambda$ is a closed conic subvariety of $T^*X$ containing
$\chr(\cM^\bullet)$. Then,
\[
\mu\eu(\cM^\bullet) =
\left[ ch_\Lambda(\sigma(\cM^\bullet))\smile \pi^*Td(TX)\right]^{2d}\ .
\]
\end{conj}
(For $\alpha$ an element of a graded object we denote by $[\alpha]^p$
the homogeneous component of $\alpha$ of degree $p$.)

We will refer to formulas such as the one in \ref{conj:ch} as Riemann-Roch
formulas and refer to the left hand side as ``non-commutative'' and the
right hand side as ``commutative''. What follows is an informal discussion
of the Riemann-Roch formula and the key ideas and observations which enter
into the statment and the proof of the main result of this paper (Theorem
\ref{thm:mainE}) from which Conjecture \ref{conj:ch} follows.

\subsubsection{}
Although Conjecture \ref{conj:ch} is stated in terms of $\cD$-modules, it
is, in fact of micro-local (i.e. local on $T^*X$) nature and is naturally
formulated in terms of modules over the ring $\cE_X$ of micro-differential
operators. Recall that $\cE_X$ is a sheaf of algebras on $T^*X$ equipped
with the canonical filtration $F_\bullet\cE_X$ by order, the symbol map
$\sigma : Gr^F_\bullet\cE_X @>>> \cO_{T^*X}$ and the canonical faithfully
flat
map $\pi^{-1}\cD_X @>>> \cE_X$ of filtered $\pi^{-1}\cO_X$-algebras.
The characteristic variety of a (coherent) $\cD_X$-module $\cM$ is the
support of the $\cE_X$-module $\pi^{-1}\cM\otimes_{\pi^{-1}\cD_X}\cE_X$.

Let $K^0_\Lambda(\cD_X)$ (respectively $K^0_\Lambda(\cE_X)$) denote the
Grothendieck group of the category of perfect complexes of $\cD_X$-modules
whose characteristic variety is contained in $\Lambda$ (respectively
perfect complexes of $\cE_X$-modules supported on $\Lambda$). Extension
of scalars gives rise to the map $K^0_\Lambda(\cD_X) @>>> K^0_\Lambda(\cE_X)$.

The microlocal Euler class which appears on the non-commutative side of
the Riemann-Roch formula is a map $K^0_\Lambda(\cD_X) @>>> H^{2d}_\Lambda
(T^*X;\bbC)$, but, in fact, by it's very definition, is the composition
of extension of scalars with the map $\mu\eu : K^0_\Lambda(\cE_X) @>>>
 H^{2d}_\Lambda (T^*X;\bbC)$.

Similarly, the commuative side of the Riemann-Roch formula depends not on
the filtered $\cD_X$-module but on the filtered $\cE_X$-module obtained
by extension of scalars.

\subsubsection{}
The presence of the filtration can be accounted for in ring theoretic terms
by the traditional device of the Rees construction. The Rees ring of the
filtered ring $(\cE_X,F_\bullet)$ is the (graded) $\bbC[t]$-algebra
$\cR\cE_X = \bigoplus_p F_p\cE_X t^p\subset\cE_X[t]$. A filtered
$\cE_X$-module is the same as a (graded) $t$-torsion-free module over
$\cR\cE_X$: $\cR$ is an exact functor defined by the formula
$\cR(\cM,F_\bullet)=\bigoplus_p F_p\cM t^p$.

The algebra $\cR\cE_X$ comes equipped with the symbol map
$\sigma : \cR\cE_X @>>> \cO_{T^*X}$ which annihilates $t$ and the canonical
isomoprhism $\cR\cE_X[t^{-1}]\isomo\cE_X[t^{-1},t]$. The symbol complex is
given by the formula
$\sigma(\cM^\bullet) = \cR\cM^\bullet\otimes_{\cR\cE}\cO_{T^*X}$
and the complex is recovered via the canonical isomorphism
$\cR\cM^\bullet [t^{-1}]\isomo\cM^\bullet [t^{-1},t]$. The natural input
for the Riemann-Roch formula is conveniently summarized as a perfect complex
of $\cR\cE_X$-modules.

The use of the Rees construction is not just a neat device to keep track of
the filtration. The appearance of $\cE_X$ and $Gr^F_\bullet\cE_X$ as,
respectively, the generic and the special fiber in a one parameter deformation
allows application of deformation theoretic metods which play a crucial role
in the proof of the Riemann-Roch formula.

\subsubsection{}
It becomes necessary, particularly with the appearance of algebras such as
$\cR\cE_X$ to use Hochschild homology, rather than appeal directly to the
duality theory for the definition of the micro-local Euler class. The map
$\mu\eu : K^0_\Lambda(\cE_X) @>>> H^{2d}_\Lambda (T^*X;\bbC)$ factors
naturally into the composition
\[
K^0_\Lambda(\cE_X) @>{\Eu^{\cE_X}_\Lambda}>>
H^0_\Lambda(T^*X;\cE_X\otimes^{\bbL}_{\cE_X\otimes\cE_X^{op}}\cE_X)
@>{\mu_{\cE}}>>
H^{2d}_\Lambda (T^*X;\bbC)
\]
where the first map is the Euler class with values in Hochschild homology
and is defined in the generality of perfect complexes of modules over sheaves
of algebras, and the second map is induced by the morphism
\[
\mu_{\cE} : \cE_X\otimes^{\bbL}_{\cE_X\otimes\cE_X^{op}}\cE_X
@>>> \bbC_{T^*X}[2d]
\]
in the derived category of sheaves on $T^*X$ called {\em the (canonical)
trace density map}. The existence of the canonical trace density map is
a non-trivial fact which has to do with special properties of the algebra
of micro-differential operators.

The above mentioned factorization of the microlocal Euler class essentially
reduces the Riemann-Roch formula to the problem of expressing the
canonical trace density map in ``commutative'' terms, i.e., loosely speaking,
factoring the composition
\[
\cR\cE_X\otimes^{\bbL}_{\cR\cE_X\otimes\cR\cE_X^{op}}\cR\cE_X
@>>>
\cE_X\otimes^{\bbL}_{\cE_X\otimes\cE_X^{op}}\cE_X[t^{-1},t]
@>{\mu_{\cE}}>> \bbC_{T^*X}[t^{-1},t][2d]
\]
through the symbol map. (This reduction is not so surprising, considering the
fact that the relevant portion of the Hochschild homology of micro-differential
operators is spanned by the Euler classes.)

Before addressing the latter, simpler, problem we make some more observations
with regard to the Riemann-Roch formula intended to motivate a gearchange that
follows.

\subsubsection{}
The added complication in the Riemann-Roch formula is the fact that
commutative side is comprised of a homogeneous component of a manifestly
non-homogeneous expression. As much as one would like to remove the
brackets in the formula it is not clear how to extend the micro-local Euler
class to a manifestly non-homogeneous characteristic class without the
reformulation of the problem in terms of Hochschild homology.

Hochschild homology is the natural recipient of the Lefshetz trace map.
The Euler class is defined as the Lefshetz trace of the identity and
posesses special symmetries not shared by the traces of other endomorphisms.
As a result, the Euler class extends to the Chern character
\[
ch^{\cE_X}_\Lambda : K^0_\Lambda(\cE_X) @>>>
H^0_\Lambda(T^*X;CC^{per}_\bullet(\cE_X)) \ .
\]

On the non-commutative side of the Riemann-Roch formula the transition from
Hochschild homology to periodic cyclic homology amounts to switching from
the integer grading to the even/odd grading, i.e. allowing not necessarily
homogeneous classes which is a welcome change. Due to the same special
homological properties of micro-differential operators the canonical trace
density map extends in a unique fashion to the morphism (in the derived
category of sheaves on $T^*X$)
\[
\tilde\mu_{\cE} : CC^{per}_\bullet(\cE_X) @>>> 
\bbC_{T^*X}[u^{-1},u]]
\]
where $u$ denotes a formal variable of degree $2$ (so that, for a complex
$A^\bullet$ multiplication by $u$ establishes the isomorphism
$A^\bullet @>>> uA^\bullet [2]$).

\subsubsection{}
The preceeding discussion was intended as a justification of the following
reformulation of the Riemann-Roch formula in terms of the (non-commutative)
Chern character: for a perfect complex $\cM^\bullet$ of $\cR\cE_X$-modules
\[
\tilde\mu_{\cE}(ch^{\cE_X}_\Lambda(\iota\cM^\bullet)) =
ch_\Lambda(\sigma\cM^\bullet)\smile\pi^*Td(TX)\ ,
\]
where $\iota\cM^\bullet = \cM^\bullet [t^{-1}]$, $\sigma\cM^\bullet =
\cM^\bullet\otimes^{\bbL}_{\cR\cE_X}\cO_{T^*X}$, and the equality takes
place in $H^{2d}_\Lambda(T^*X;\bbC_{T^*X}[t^{-1},t][u^{-1},u]])\isomo
H^{even}_\Lambda(T^*X;\bbC_{T^*X}[t^{-1},t])$.

\subsubsection{}
The Chern character which appears on the commutative side of the
Riemann-Roch formula may be formulated in terms of periodic cyclic
homology of $\cO_{T^*X}$ as the composition
\[
K^0_\Lambda(\cO_{T^*X}) @>{ch^{\cO_{T^*X}}_\Lambda}>>
H^0_\Lambda(T^*X;CC^{per}(\cO_{T^*X}) @>{\tilde\mu_{\cO}}>>
H^0_\Lambda(T^*X;\bbC_{T^*X}[u^{-1},u]])
\]
where the second map is induced by the Hochschild-Kostant-Rosenberg-Connes
map $\tilde\mu_{\cO} : CC^{per}(\cO_{T^*X}) @>>>
\Omega^\bullet_{T^*X}[u^{-1},u]]$ and the de Rham isomorphism.

With these notations the Riemann-Roch formula reads
\[
\tilde\mu_{\cE}(ch^{\cE_X}_\Lambda(\iota\cM^\bullet)) =
\tilde\mu_{\cO}(ch^{\cO_{T^*X}}_\Lambda(\sigma\cM^\bullet))
\smile\pi^*Td(TX)\ .
\]

Furthermore, the naturality of the Chern character implies that
$ch^{\cE_X}_\Lambda(\iota\cM^\bullet) =
\iota(ch^{\cE_X}_\Lambda(\cM^\bullet))$ (respectively
$ch^{\cO_{T^*X}}_\Lambda(\sigma\cM^\bullet) =
\sigma(ch^{\cO_{T^*X}}_\Lambda(\cM^\bullet))$) where
$\iota : CC^{per}_\bullet(\cR\cE_X) @>>>  CC^{per}_\bullet(\cE_X)[t^{-1},t]$
(respectively $\sigma : CC^{per}_\bullet(\cR\cE_X) @>>> 
CC^{per}_\bullet(\cO_{T^*X})$) is induced by
$\iota : \cR\cE_X @>>> \cE_X[t^{-1},t]$ (respectively
$\sigma : \cR\cE_X @>>> \cO_{T^*X}$). The Riemann-Roch formula now takes
the shape
\[
\tilde\mu_{\cE}(\iota(ch^{\cE_X}_\Lambda(\cM^\bullet))) =
\mu_{\cO}(\sigma(ch^{\cO_{T^*X}}_\Lambda(\cM^\bullet)))\smile\pi^*Td(TX)
\]
and, in the present formulation is valid for any cycle in $CC^{per}_\bullet
(\cR\cE_X)$ (because the image of the canonical trace density map is, in
essence, spanned by Chern characters). The above formula is equivalent to
the commutativity of the diagram
\[
\begin{CD}
CC^{per}_\bullet(\cR\cE_X) @>{\sigma}>> CC^{per}_\bullet(\cO_{T^*X}) \\
@V{\iota}VV @VV{\tilde\mu_{\cO}\smile\pi^*Td(TX)}V \\
CC^{per}_\bullet(\cE_X)[t^{-1},t] @>{\tilde\mu_{\cE}}>>
\bbC_{T^*X}[t^{-1},t][u^{-1},u]]
\end{CD}
\]
in the derived category of sheaves on $T^*X$.

\subsubsection{}
The non-triviality of the issue of defining the canonical trace density map
has to do with the fact that the canonical morphism $\bbC_{T^*X}[2d] @>>>
\cE_X\otimes^{\bbL}_{\cE_X\otimes\cE_X^{op}}\cE_X$ furnished by the duality
theory is not an isomorphism (in the derived category) and that can be
traced to the fact that the canonical map $Gr^F_\bullet\cE_X @>>>
\cO_{T^*X}$ is not an isomorphism.

The trace density map is a retraction of the morphism $\bbC_{T^*X}[2d] @>>>
\cE_X\otimes^{\bbL}_{\cE_X\otimes\cE_X^{op}}\cE_X$ and is arranged for by
embedding $\cE_X$ into a larger algebra. The choice more natural in the context
of micro-local analysis is the algebra of micro-local operators
$\cE_X^{\bbR}$ which bears very little resemblance to $\cO_{T^*X}$.

The same may be achieved by embedding $\cR\cE_X$ into a formal one-parameter
deformation (with $t$ being the parameter)
$\bbA^t_{T^*X}$ of $\cO_{T^*X}$ which is uniquely determined by the fact that
it contains $\cR\cE_X$. The Riemann-Roch formula in its diagramatic
reformulation, follows from the analogous statement about $\bbA^t_{T^*X}$:
\[
\begin{CD}
CC^{per}_\bullet(\bbA^t_{T^*X}) @>{\sigma}>{t\mapsto 0}>
CC^{per}_\bullet(\cO_{T^*X}) \\
@V{\iota}VV @VV{\tilde\mu_{\cO}\smile\pi^*Td(TX)}V \\
CC^{per}_\bullet(\bbA^t_{T^*X})[t^{-1}] @>{\tilde\mu^{\bbA}}>>
\bbC_{T^*X}[t^{-1},t]][u^{-1},u]]
\end{CD}
\]
commutes in the derived category of sheaves on $T^*X$. The latter fact is
a particular case of a general Riemann-Roch type theorem conserning
symplectic deformations of structure sheaves and whose proof relies on
the techniques of formal geometry and invariant theory and constitutes
the sequel to the present paper.

\subsection{Acknowledgements}
The authors would like to thank J.-L.Brylinski and one another for
inspiring discussions.

\section{Characteristic classes and trace maps}\label{section:outline}
In this section we introduce the ingredients which go into the statement of
the main result of this paper (Theorem \ref{thm:main}) and show how
Conjecture \ref{conj:ch} follows from it. Proofs are postponed until later
sections as indicated. 

\subsection{Characteristic classes of perfect complexes}
The following theorem (Theorem \ref{thm-intro:char-cl}) summarizes the relevant
aspects of the requisite characteristic classes for perfect complexes.
The details of the construction will appear in a separate publication.

Let $X$ be a topological space, $Z$ a closed subset of $X$.
Let $\cA$ denote a sheaf of algebras on $X$ such that there is a global
section $1\in\Gamma(X;\cA)$ which restricts to $1_{\cA_x}$.

We denote by $K^i_Z(\cA)$ the $i$-th $K$-group of the category of perfect
complexes of $\cA$-modules which are acyclic on the complement of $Z$ and
refer the reader to Section \ref{section:CC} for other notations.

\begin{thm}\label{thm-intro:char-cl}
There exists the Chern character
$ch^{\cA}_{Z,i} : K^i_Z(\cA)\to H^{-i}_Z(X;CC^-_\bullet(\cA))$ and the Euler class
$\Eu^{\cA}_{Z,i} : K^i_Z(\cA)\to H^{-i}_Z(X;C_\bullet(\cA))$ natural in $X$, $Z$,
and $\cA$ such that
\begin{itemize}
\item the composition $K^i_Z(\cA) @>{ch^{\cA}_{Z,i}}>> H^{-i}_Z(X;CC^-_\bullet(\cA))
@>>> H^{-i}_Z(X;C_\bullet(\cA))$ coincides with $\Eu^{\cA}_{Z,i}$;
\item for a perfect complex $\cF^\bullet$ of $\cA$-modules supported on $Z$
the class $\Eu^{\cA}_{Z,0}(\cF^\bullet)\in H^0_Z(X;C_\bullet(\cA))$ coincides with
the compostion
\begin{multline*}
k @>{1\mapsto\id}>> \R\shHom^\bullet_{\cA}(\cF^\bullet,\cF^\bullet)
@<{\isomo}<< (\R\shHom^\bullet_{\cA}(\cF^\bullet,\cA)\otimes_k\cF^\bullet)
\otimes^{\bbL}_{\cA\otimes_k\cA^{op}}\cA \\ @>{ev\otimes\id}>>
\cA\otimes^{\bbL}_{\cA\otimes_k\cA^{op}}\cA\ .
\end{multline*} 
\end{itemize}
\end{thm}

In what follows we will only consider the components of the Euler class and
the Chern character defined on $K^0$ and write $\Eu^{\cA}_Z$ (respectively
$ch^{\cA}_Z$) instead of $\Eu^{\cA}_{Z,0}$ (respectively
$ch^{\cA}_{Z,0}$). Nor shall we make notational distinction between $ch^{\cA}_Z$
and its image under the natural map $H^\bullet_Z(X;CC^-_\bullet(\cA))\to
H^\bullet_Z(X;CC^{per}_\bullet(\cA))$.

\subsection{The microlocal Euler class}
For the reader's convenience we give the definition of the microlocal
Euler class. 

Suppose that $\Lambda$ is a closed subvariety of $T^*X$ and
$\cN^\bullet$ is a perfect complex of $\cE^{\bbR}_X$-modules acyclic on the
complement of $\Lambda$. Consider the morphism in the derived category
(the Lefshetz trace map) given by the composition 
\begin{eqnarray*}
\R\shHom_{\cE^{\bbR}_X}(\cN^\bullet, \cN^\bullet) & \isomo &
\cN^\bullet\underline{\boxtimes}\underline{D}\cN^\bullet
\otimes^{\bbL}_{\cE^{\bbR}_{X\times X}}\cC^{\bbR}_{\Delta\vert X\times X} \\
& @<{\isomo}<< & \R\Gamma_\Lambda\left(
\cN^\bullet\underline{\boxtimes}\underline{D}\cN^\bullet
\otimes^{\bbL}_{\cE^{\bbR}_{X\times X}}\cC^{\bbR}_{\Delta\vert X\times X}
\right) \\
& @>>> & \R\Gamma_\Lambda\left(\cC^{(2d)\bbR}_{\Delta\vert X\times X}
\otimes^{\bbL}_{\cE^{\bbR}_{X\times X}}\cC^{\bbR}_{\Delta\vert X\times X}
\right) \\
& \isomo &  \R\Gamma_\Lambda\left(\bbC_{T^*X}[2d]\right)
\end{eqnarray*}
where
\[
\underline{D}\cN^\bullet\overset{def}{=}
\R\shHom_{\cE^{\bbR}_X}(\cN^\bullet,\cE^{\bbR}_X\otimes_{\pi^{-1}\cO_X}
\pi^{-1}\Omega^d_X[d])\ ,
\]
$\pi : T^*X @>>> X$ is the canonical projection, and the last isomorphism
is as in ... . 

The microlocal Euler class $\mu\eu(\cN^\bullet)$ is defined
as the Lefshetz trace of the identity, i.e. as the composite
\[
\mu\eu(\cN^\bullet) : \bbC_{T^*X} @>{1\mapsto\id}>>
\R\shHom_{\cE^{\bbR}_X}(\cN^\bullet, \cN^\bullet) @>>> 
 \R\Gamma_\Lambda\left(\bbC_{T^*X}[2d]\right)\ .
\]

For a perfect complex $\cN^\bullet$ of $\cE_X$-modules the microlocal Euler
class is defined by
\[
\mu\eu(\cN^\bullet)\overset{def}{=}\mu\eu(\cN^\bullet
\otimes_{\cE_X}\cE^{\bbR}_X)\ .
\]

Similarly, the microlocal Euler class of a complex $\cN^\bullet$ of
$\cD_X$-modules with $\chr\cN^\bullet\subseteq\Lambda$ is defined as
$\mu\eu(\cN^\bullet\otimes_{\pi^{-1}\cD_X}\cE^{\bbR}_X)$.

The microlocal Euler class gives rise to the map
\[
\mu\eu : K^0_\Lambda(\cE_X) @>>> H^{2d}_\Lambda(T^*X;\bbC)\ .
\]

\subsection{The canoncial trace density map}
Let $\Delta\subset X\times X$ denote the diagonal.
Let $a: T^*X @>>> T^*X$ denote the antipodal map ($(x,\xi)\mapsto (x,-\xi)$).
There is a natural isomorphism of sheaves of algebras $a^{-1}\cE^{(\bbR)}_X
\isomo\left(\cE^{(\bbR)}_X\right)^{op}$. Put
$\left(\cE^{(\bbR)}_{X\times X}\right)^a
= (\id\times a)^{-1} \cE^{(\bbR)}_{X\times X}$.
Then, there is a natural map 
$\cE^{(\bbR)}_X\boxtimes\left(\cE^{(\bbR)}_X\right)^{op} @>>>
\left(\cE^{(\bbR)}_{X\times X}\right)^a$, and the canonical structure of an
$\left(\cE^{(\bbR)}_{X\times X}\right)^a$-module on $\cE^{(\bbR)}_X$ which
extends the canoncial structure of an $\cE^{(\bbR)}_X\boxtimes
\left(\cE^{(\bbR)}_X\right)^{op}$-module.

By a theorem of Kashiwara the canonical map
\[
\bbC_{T^*X} @>>> \R\shHom_{\cE_{X\times X}^{\bbR}}(
\cC_{\Delta\vert X\times X}^{\bbR},
\cC_{\Delta\vert X\times X}^{\bbR})
\]
is an isomorphism in the derived category.
There are canonical isomorphisms
\begin{eqnarray*}
\bbC_{T^*X}[2d] & \isomo & 
\R\shHom_{\cE_{X\times X}^{\bbR}}(
\cC_{\Delta\vert X\times X}^{\bbR},
\cC_{\Delta\vert X\times X}^{\bbR})[2d] \\
& \isomo & \R\shHom_{\cE_{X\times X}^{\bbR}}
(\cC_{\Delta\vert X\times X}^{\bbR},\cE_{X\times X}^{\bbR})[2d]
\otimes^{\bbL}_{\cE_{X\times X}^{\bbR}}
\cC_{\Delta\vert X\times X}^{\bbR} \\
& \isomo &
\cC^{(2d)\bbR}_{\Delta\vert X\times X}\otimes^{\bbL}_{\cE_{X\times X}^{\bbR}}
\cC_{\Delta\vert X\times X}^{\bbR} \\
& \isomo & \cE^{\bbR}_X\otimes^{\bbL}_{\left(\cE^{\bbR}_{X\times X}\right)^a}
\cE^{\bbR}_X\ .
\end{eqnarray*}

By abuse of notation we will denote by $\mu_{\cE}$ and refer to as {\em the canonical
trace density map} the composition
\[
\cE_X\otimes^{\bbL}_{\left(\cE_{X\times X}\right)^a}\cE_X @>>>
\cE^{\bbR}_X\otimes^{\bbL}_{\left(\cE^{\bbR}_{X\times X}\right)^a}
\cE^{\bbR}_X @>{\isomo}>> \bbC_{T^*X}[2d]
\]
and the composition of the latter with the canonical map
\[
\cE_X\otimes^{\bbL}_{\cE_X\boxtimes\cE_X^{op}} \cE_X @>>>
\cE_X\otimes^{\bbL}_{\left(\cE_{X\times X}\right)^a}\cE_X\ .
\]

\subsection{Comparison of the Euler classes}
The Euler class of Theorem \ref{thm-intro:char-cl} combined with the canonical 
trace density map gives rise to the map
\begin{equation}\label{map:tr-Eu}
K^0_\Lambda(\cE_X) @>{\Eu^{\cE_X}}>>
H^0_\Lambda(T^*X;\cE_X\otimes^{\bbL}_{\cE_X\boxtimes\cE_X^{op}}\cE_X)
@>{\mu_{\cE}}>> H^{2d}_\Lambda(T^*X; \bbC)\ .
\end{equation}

\begin{prop}\label{prop:mueu-is-Eu}
The composition \eqref{map:tr-Eu} coincides with the micro-local Euler class.
In other words, if $\cN^\bullet$ is a perfect complex of $\cE_X$-modules and
$\Lambda$ is a closed subvariety of $T^*X$ containing $\supp\cN^\bullet$,
then $\mu\eu(\cN^\bullet) = \mu_{\cE}(\Eu^{\cE_X}_\Lambda(\cN^\bullet))$.
\end{prop}
\begin{pf}
The Lefshetz trace map
\[
\R\shHom_{\cE_X}(\cN^\bullet, \cN^\bullet) @>>>
\cC^{(2d)}_{\Delta\vert X\times X}\otimes^{\bbL}_{\cE_{X\times X}}
\cC_{\Delta\vert X\times X}
\]
is compatible with the Lefshetz trace map
\[
\R\shHom_{\cE_X}(\cN^\bullet, \cN^\bullet) @>>>
\cE_X\otimes^{\bbL}_{\cE_X\otimes\cE_X^{op}}\cE_X
@>>>
\cE_X\otimes^{\bbL}_{\left(\cE_{X\times X}\right)^a}\cE_X
\]
and the isomorphisms
\[
\cC^{(2d)\bbR}_{\Delta\vert X\times X}\otimes^{\bbL}_{\cE^{\bbR}_{X\times X}}
\cC^{\bbR}_{\Delta\vert X\times X}\isomo
\cE^{\bbR}_X\otimes^{\bbL}_{\left(\cE^{\bbR}_{X\times X}\right)^a}\cE^{\bbR}_X
\isomo \bbC_{T^*X}[2d]\ .
\]
\end{pf}

\subsection{Complexes of Hochschild chains}
With the view onto the passage to cyclic homology we describe the canonical
trace density map in terms of the standard complexes of Hochschild chains.

The object $\cE_X\otimes^{\bbL}_{\cE_X\boxtimes\cE_X^{op}}\cE_X$, regarded
as an object of the derived category of sheaves on $T^*X$, is represented by
the standard complex of Hochschild chains $C_\bullet(\cE_X)$
(see \ref{subsection:C-CC}). Put
\begin{eqnarray*}
\left(\cE^{(\bbR)}_X\right)^{\frak{e}} &= &
\delta_{T^*X}^{-1}(\id\times a)^{-1}
\cE^{(\bbR)}_{X\times X} \\
\widehat{C}_p(\cE^{(\bbR)}_X) & = & \delta_{T^*X}^{-1}
\cE^{(\bbR)}_{X^{\times p+1}}\ ,
\end{eqnarray*}

The Hochschild differential extends to the map
$b: \widehat{C}_p(\cE^{(\bbR)}_X) @>>> \widehat{C}_{p-1}(\cE^{(\bbR)}_X)$.
The complex
$\widehat{C}_\bullet(\cE^{(\bbR)}_X)$ represents
$\cE^{(\bbR)}_X\otimes^{\bbL}_{\left(\cE^{(\bbR)}_X\right)^{\frak{e}}}
\cE^{(\bbR)}_X$ in the derived category of sheaves on $T^*X$. There is a
canoncial map $C_\bullet(\cE^{(\bbR)}_X) @>>>
\widehat C_\bullet(\cE^{(\bbR)}_X)$.

Let $\Phi^{\cE}$ denote the global section of
$H^{-2d}\widehat{C}_\bullet(\cE^{\bbR}_X)$
which corresponds to the global section $1$ under the
isomorphism induced by $\mu_{\cE}$.

\begin{lemma}
Let $x_1,\ldots,x_d$ be a local coordinate system on $X$ and let
$\partial_i =\displaystyle\frac{\partial}{\partial x_i}$. Then,
the locally defined Hochschild chain
$Alt(1\otimes x_1\otimes\cdots\otimes x_d
\otimes\partial_1\otimes\cdots\otimes\partial_d)\in C_{2d}(\cE_X)$
is a cycle which represents (the restriction to the domain of the coordinate
system of)
$\Phi^{\cE}\in\Gamma (T^*X;H^{-2d}C_\bullet(\cE^{\bbR}_X))$.
\end{lemma}

\subsection{Cyclic homology of microdifferential operators}\label{subsec:CCE}
Our immediate goal is to recast the calculation of the micro-local Euler
class in terms of the negative (and, ultimately, the periodic) cyclic homology.
The passage to cyclic homology is motivated by two developments. 

On the one hand, the Euler class (which takes values in Hochschild homology)
factors canonically through the Chern character which takes values in
negative cyclic homology; thus, the Chern character should be a (significant)
refinement of the Euler class. On the other hand, the periodic cyclic homology
is invariant under formal deformations of algebras and this feature will,
ultimately, allow us to evaluate the ``non-commutative'' Chern character
explicitly in ``commutative'' terms.

The cyclic differential $B$ extends to the map
$B: \widehat{C}_p(\cE^{(\bbR)}_X) @>>> \widehat{C}_{p+1}(\cE^{(\bbR)}_X)$
and the cyclic complexes $\widehat{CC}^-_\bullet(\cE^{(\bbR)})$ and
$\widehat{CC}^{per}_\bullet(\cE^{(\bbR)}_X)$ are defined as in 
\ref{subsection:C-CC} starting from $\widehat C_\bullet(\cE^{(\bbR)}_X)$.
The analogs of \eqref{inclusions} and \eqref{ses:CC-C} hold and there are
canonical maps
${CC}^-_\bullet(\cE^{(\bbR)}) @>>> \widehat{CC}^-_\bullet(\cE^{(\bbR)})$
and
${CC}^{per}_\bullet(\cE^{(\bbR)}_X) @>>>
\widehat{CC}^{per}_\bullet(\cE^{(\bbR)}_X)$.

The complexes $(\widehat{CC}^-_\bullet(\cE^{(\bbR)}), b+uB)$ and
$(\widehat C_\bullet(\cE^{(\bbR)}_X)[[u]], b)$ have canonical filtrations
by powers of $u$ which will be denoted by $F_\bullet$. Note that the associated
graded complexes $Gr^F_\bullet\widehat{CC}^-_\bullet(\cE^{(\bbR)})$ and
$Gr^F_\bullet\widehat C_\bullet(\cE^{(\bbR)}_X)[[u]]$ are identical.

\begin{lemma}\label{lemma:CC-E}
There exists an isomorphism of filtered complexes
\[
(\widehat{CC}^{-}_\bullet(\cE^{\bbR}_X),b+uB, F_\bullet)\isomo
(\widehat C_\bullet(\cE^{\bbR}_X)[[u]],b,F_\bullet)
\]
which induces the identity map on the associated graded complexes and any two
such are (filtered) homotopic.
\end{lemma}
\begin{pf}
Consider the spectral sequence of the filtered complex \\
$(\Hom^\bullet(\widehat{CC}^{-}_\bullet(\cE^{\bbR}_X),
\widehat C_\bullet(\cE^{\bbR}_X)[[u]]),F_\bullet)$ and use the isomorphism
(in the derived category)
$\widehat C_\bullet(\cE^{\bbR}_X)\isomo\bbC_{T^*X}[2d]$.
\end{pf}
\begin{cor}\label{cor:CC-E}
There exists an isomorphism of filtered complexes
\[
(\widehat{CC}^{per}_\bullet(\cE^{\bbR}_X),b+uB,F_\bullet)
\isomo
(\widehat C_\bullet(\cE^{\bbR}_X)[u^{-1},u]],b,F_\bullet)
\]
which induces the identity map on the associated graded objects and any two
such are homotopic.
\end{cor}

\subsection{The trace density on the periodic cyclic complex}
\label{subsec:trace-density}
The isomorphism (in the derived category)
$\widehat{C}_\bullet(\cE^{\bbR}_X)
\isomo\bbC_{T^*X}[2d]$ is represented, using the standard truncation functors,
by the following diagram of quasi-isomorphisms of complexes:
\begin{multline}\label{map:muE-trunk}
\widehat{C}_\bullet(\cE^{\bbR}_X)
@<<< 
\tau^{\leq -2d}\widehat{C}_\bullet(\cE^{\bbR}_X) @>>> \\
@>>>
\tau^{\geq -2d}\tau^{\leq -2d}\widehat{C}_\bullet(\cE^{\bbR}_X)
\isomo
\bbC_{T^*X}[2d]\ .
\end{multline}

The diagram \eqref{map:muE-trunk} gives rise to the diagram of complexes of
$\bbC[u^{-1},u]]$-modules
\begin{multline*}
\widehat{C}_\bullet(\cE^{\bbR}_X)[u^{-1},u]]
@<<< 
\left(\tau^{\leq -2d}\widehat{C}_\bullet(\cE^{\bbR}_X)\right)[u^{-1},u]] @>>> \\
@>>>
\left(\tau^{\geq -2d}\tau^{\leq -2d}\widehat{C}_\bullet(\cE^{\bbR}_X)\right)
[u^{-1},u]]
\isomo
\bbC_{T^*X}[2d][u^{-1},u]]
\end{multline*}
which represents a morphism
\begin{equation}\label{map:muE-trunk-u}
\widehat{C}_\bullet(\cE^{\bbR}_X)[u^{-1},u]] @>>>
\bbC_{T^*X}[2d][u^{-1},u]]
\end{equation}
in the derived category.

We will denote by $\tilde\mu_{\cE}$ and refer to as the
canonical trace density map the composition
\begin{multline*}
{CC}^{per}_\bullet(\cE_X)
@>>>
{CC}^{per}_\bullet(\cE^{\bbR}_X)
@>>> \\
\widehat{CC}^{per}_\bullet(\cE^{\bbR}_X)
@>>> 
(\widehat C_\bullet(\cE^{\bbR}_X)[u^{-1},u]]
@>>>
\bbC_{T^*X}[2d][u^{-1},u]]
\end{multline*}
where the last two maps are furnished, respectively, by Corollary
\ref{cor:CC-E} and \eqref{map:muE-trunk-u}. 

Note that, by construction, the map $\tilde\mu_{\cE}$ induces the map
$\mu_{\cE}$ on the associated (to the filtrations by powers of $u$) graded
objects.

\begin{lemma}\label{cor:Eu-comp-ch}
For a perfect complex $\cN^\bullet$ of $\cE_X$-modules
and $\Lambda$ a closed subvariety of $T^*X$ containing $\supp\cN^\bullet$,
$\mu_{\cE}(\Eu^{\cE_X}_\Lambda(\cN^\bullet)) = \left[\tilde\mu_{\cE}
(ch^{\cE_X}_\Lambda(\cN^\bullet))\right]^{2d}$.
\end{lemma}
\begin{pf}
Consider the commutative diagram
\[
\begin{CD}
H^0_\Lambda(T^*X;CC^-_\bullet(\cE_X)) @>{\tilde\mu_{\cE}}>>
H^0_\Lambda(T^*X;\bbC_{T^*X}[2d][[u]]) \\
@VVV @VV{u\mapsto 0}V \\
H^0_\Lambda(T^*X;C_\bullet(\cE_X)) @>{\mu_{\cE}}>>
H^{2d}_\Lambda(T^*X;\bbC)
\end{CD}
\]
\end{pf}

\section{The Schapira--Schneiders conjecture}

In this section we arrive at the statement of the main result of this paper
which implies the Schapira--Schneiders conjecture and explain how it follows
from the the analogous theorem for formal symplectic deformations. The proof
of the latter is presented in the sequel to this paper.

\subsection{The Rees construction}
The sheaf of algebras $\cE_X$ carries a natural filtration $F_\bullet\cE_X$
by order. The Rees ring $\cR\cE_X$ is the (graded) algebra flat over
$\bbC[\hbar]$ defined by
\[
\cR\cE_X\overset{def}{=}\bigoplus_p F_p\cE_X\cdot t^p\ .
\]
The Rees module $\cR\cN$ (which is a graded module over $\cR\cE_X$) associated
to a filtered module $(\cN,F_\bullet)$ is defined similarly.
Thus defined, the Rees construction extends to an exact functor $\cR$ form
the (exact) category of filtered $\cE_X$-modules to the (Abelian) category
of graded $\cR\cE_X$-modules. The functor $\cR$ is an embedding with the
essential image consisting of the subcategory of $t$-torsion free modules.
The induced functor between the respective derived categories is an
equivalence.

The filtered complex $(\cN^\bullet,F_\bullet)$ is {\em good} if and only if
the complex of $\cR\cE_X$-modules $\cR\cN^\bullet$ is perfect. The Rees
construction restricts to an equivalence between the derived categories of
good filtered complexes of $\cE_X$-modules and perfect complexes of graded
$\cR\cE_X$-modules. In particular it induces an isomorphism of respective
Grothendieck groups.

Setting $t = 0$ one obtains the map
\[
\sigma : \cR\cE_X @>>> Gr^F_\bullet\cE_X @>>> \cO_{T^*X}\ .
\]
Note that
\[
\sigma(\cN)\overset{def}{=}Gr^F_\bullet\cN\otimes_{Gr^F_\bullet\cE_X}
\cO_{T^*X} = \cR\cN\otimes_{\cR\cE_X}\cO_{T^*X}\ .
\]

Consider a filtered complex $(\cN^\bullet,F_\bullet)$ with $\cR\cN^\bullet$
perfect over $\cR\cE_X$ with support contained in $\Lambda\subset T^*X$. Then,
\[
ch^{\cO_{T^*X}}_\Lambda(\sigma(\cN^\bullet)) = 
\sigma\left(ch^{\cR\cE_X}_\Lambda(\cR\cN^\bullet)\right)\ .
\]

The canonical homomorphism of algebras
\[
\iota :\cR\cE_X @>>> \cR\cE_X[t^{-1}] @>{\isomo}>> \cE_X[t^{-1},t]
\]
induces the natural isomorphism $\cR\cN^\bullet[t^{-1}]\isomo\cN^\bullet
\otimes_{\cE_X}\cE_X[t^{-1},t]$. Therefore,
\[
ch^{\cE_X}_\Lambda(\cN^\bullet) = 
\iota\left(ch^{\cR\cE_X}_\Lambda(\cR\cN^\bullet)\right) \ .
\]
In particular the right hand side is independent of $t$.

\subsection{Cyclic homology of the structure sheaf}
We recall, briefly, the well known calculation of the cyclic homology of the
structure sheaf $\cO_M$ of a complex manifold $M$ of
Hochschild-Kostant-Rosenberg-Connes.

Let $\widehat C_p(\cO_M)$ denote the completion of $\cO_{M^{\times p+1}}$
along the diagonal. The Hochschild differential extends to the map
$b: \widehat C_p(\cO_M) @>>> \widehat C_{p-1}(\cO_M)$ and the resulting complex
$\widehat C_\bullet(\cO_M)$ represents $\cO_M\otimes^{\bbL}_{\cO_{M\times M}}\cO_M$
in the derived category. There is a canonical map $C_\bullet(\cO_M) @>>>
\widehat C_\bullet(\cO_M)$.

The assignment $f_0\otimes\dots\otimes f_p\mapsto\displaystyle\frac1{p!}
f_0df_1\wedge\dots df_p$ extends to a of complexes
\[
\mu_{\cO} : \widehat C_\bullet(\cO_M) @>>> \bigoplus_p\Omega^p_M[p]
\]
which, according to a theorem of Hochschild-Kostant-Rosenberg, is
a quasi-isomorphism.

The cyclic differential $B$ extends to the map
$B : \widehat C_p(\cO_M) @>>> \widehat C_{p+1}(\cO_M)$ so that the square
\[
\begin{CD}
\widehat C_p(\cO_M) @>>> \Omega^p_M \\
@V{B}VV @VV{d}V \\
\widehat C_{p+1}(\cO_M) @>>> \Omega^{p+1}_M
\end{CD}
\]
commutes. The cyclic complex $\widehat{CC}^{per}_\bullet(\cO_X)$ is defined as
in \ref{subsection:C-CC}. Thus, the Hochschild-Kostant-Rosenberg map induces
the map of complexes
\[
\tilde\mu_{\cO} : \widehat{CC}^{per}_\bullet(\cO_X) @>>>
\Omega^\bullet_M[u^{-1},u]]
\]
which, according to A.~Connes, is a quasi-isomorphism.

The natural inclusion $\bbC_M[u^{-1},u]] @>>> CC^{per}(\bbC_M)$ is easily seen
to be a quasi-isomorphism and the composition
\begin{multline*}
\bbC_M[u^{-1},u]] @>>>
CC^{per}(\bbC_M) @>>> \\
CC^{per}_\bullet(\cO_X) @>>> 
\widehat{CC}^{per}_\bullet(\cO_X) @>>> 
\Omega^\bullet_M[u^{-1},u]]
\end{multline*}
coincides with the canoncal map $\bbC_M @>>> \Omega^\bullet_M$.

By abuse of notation we will denote by $\tilde\mu_{\cO}$ the morphism in the
derived category represented by
\begin{equation}\label{map:muO}
\tilde\mu_{\cO} : CC^{per}_\bullet(\cO_X) @>>> \Omega^\bullet_M[u^{-1},u]]
@<{\isomo}<< \bbC_M[u^{-1},u]]
\end{equation}
which, by the discussion above, is a one-sided inverse to the map induced by
the inclusion $\bbC_M @>>> \cO_M$.

The Chern character with supports (as in \eqref{map:ch-with-supp}) for perfect
complexes of $\cO_X$-modules supported on a closed subvariety $Z$ of a
complex manifold $M$ is defined
as the composition
\[
K^0_Z(\cO_X) @>{ch^{\cO_X}_Z}>> H^0_Z(X;CC^{per}_\bullet(\cO_X))
@>{\tilde\mu_{\cO}}>> \bigoplus_{p}H_Z^{2p}(X;\bbC)\ .
\]

\subsection{The Riemann-Roch formula}
In view of Lemma \ref{cor:Eu-comp-ch} and Proposition
\ref{prop:mueu-is-Eu}, to prove Conjecture \ref{conj:ch} it is sufficient
to show that
\begin{multline}\label{formula:RR-for-E}
\tilde\mu_{\cE}\left(\iota\left(ch^{\cR\cE_X}_\Lambda(\cR\cN^\bullet)\right)
\right) =
\tilde\mu_{\cO}\left(\sigma\left(ch^{\cR\cE_X}_\Lambda(\cR\cN^\bullet)\right)
\right)\smile\pi^*Td(TX)\ .
\end{multline}

The above reformulation of Conjecture \ref{conj:ch} in terms of the
Chern character with values in periodic cyclic homology is crucial
since, as it turns out, the formula \eqref{formula:RR-for-E} holds
with $ch^{\cR\cE_X}_\Lambda(\cR\cN^\bullet)$ replaced by an arbitrary
cycle in the periodic cyclic complex. The latter fact constitutes
the main result of the paper and is formulated in the following theorem.

\begin{thm}\label{thm:mainE}
The diagram in the derived category of sheaves on $T^*X$
\[
\begin{CD}
CC^{per}_\bullet(\cR\cE_X) @>{\sigma}>> CC^{per}_\bullet(\cO_{T^*X}) \\
@V{\iota}VV  @VV{\tilde\mu_{\cO}\smile \pi^*Td(TX)}V \\
CC^{per}_\bullet(\cE_X)[t^{-1},t] @>{\tilde\mu_{\cE}}>> \bbC_{T^*X}[t^{-1},t]
[u^{-1},u]]
\end{CD}
\]
is commutative.
\end{thm}
\begin{pf}
The statement follows from Theorem \ref{thm:main} applied to the case
$M=T^*X$ and $\bbA^t_{T^*X}$ furnished by Proposition \ref{prop:DQcotan}
(with $\theta = \pi^*c_1(TX)/2$), Corollary \ref{cor:EtoA}, and the equality
$\widehat A(TM)\smile e^\theta = \pi^*Td(TX)$.
\end{pf}

\section{Deformation quantization}
\subsection{Weyl quantization}
Consider the vector space $\bbC^{2d}$ with coordinates $x_1,\ldots,x_d,\xi_1,
\ldots,\xi_d$ as a symplectic manifold equipped with the standard symplectic
form $dx_1\wedge d\xi_1\wedge\dots\wedge dx_d\wedge d\xi_d$.
The sheaf of algebras $\bbA^t_{\bbC^{2d}}$ is defined as follows.

For $U$ an open subset of $\bbC^{2d}$ the $\bbC[[t ]]$-module underlying
$\bbA^t_{\bbC^{2d}}(U)$ is $\cO_{\bbC^{2d}}(U)[[t]]$ and
the product on $\bbA^t_{\bbC^{2d}}(U)$ is given by the standard
Moyal--Weyl product
\[
(f\ast g)(\underline x,\underline\xi) = \\
exp\left(\frac{t}{2}\sum_{i=1}^d\left(
\frac{\partial\ }{\partial\xi_i}\frac{\partial\ }{\partial y_i}-
\frac{\partial~}{\partial\eta_i}\frac{\partial~}{\partial x_i}\right)\right)
f(\underline x,\underline\xi)g(\underline y,\underline\eta)
\vert_{\overset{\underline x=\underline y}{\underline\xi =\underline\eta}}
\]
where $\underline x = (x_1,\ldots,x_d),\ \underline\xi = (\xi_1,\ldots,\xi_d),
\ \underline y = (y_1,\ldots,y_d),\ \underline\eta = (\eta_1,\ldots,\eta_d)$.

\subsection{Symplectic deformation quantization}
Let $(M,\omega)$ denote a symplectic manifold of dimension $\dim_\bbC M = 2d$.
For purposes of this paper  a {\em (symplectic) deformation quantization} of
$M$ is a formal one parameter deformation of the structure sheaf $\cO_M$,
i.e. a sheaf of algebras $\bbA^t_M$ flat over $\bbC [[t]]$ equipped with an
isomorphism of algebras $\bbA^t_M\otimes_{\bbC [[t]]}\bbC @>>> \cO_M$ and
having the following properties:
\begin{itemize}
\item it is locally isomorphic to the standard deformation
$\bbA^t_{\bbC^{2d}}$ of $\bbC^{2d}$ (see below);
\item the Poisson bracket on $\cO_M$ defined by
$\displaystyle{
\lbrace f, g\rbrace = \frac{1}{t}[\tilde f,\tilde g ]
+ t\cdot \bbA^t_M}$
where $f$ and $g$ are two local sections of $\cO_M$ and $\tilde f$, $\tilde g$
are their respective lifts $\bbA^t_M$, coincides with the one induced
by the symplectic structure.
\end{itemize}
 
To a symplectic deformation quantization $\bbA^t_M$ one associates
a characteristic class $\theta\in H^2(M;\displaystyle\frac{1}{t}\bbC[[t]])$
(see \cite{NT3} and \ref{subsection:Fedosoff} below) with the property that
the coefficient of $\displaystyle\frac{1}{t}$ is the class of the symplectic form.

In what follows we will use the canonical maps $\sigma :\bbA^t_M\to\cO_M$
and $\iota:\bbA^t_M\to\bbA^t_M[t^{-1}]$.

\subsection{Review of Fedosov connections}\label{subsection:Fedosoff}
We recall briefly the definition of the characteristic class $\theta$ in terms of
the Fedosov connections.

Let $W = \bbC[[\hat{z}_1,\dots,\hat{z}_d,\hat{\xi}_1,\dots,\hat{\xi}_d,t]]$ equipped
with the Weyl product denoted $\ast$ while the ``commutative'' product will not be
indicated by any special notations. Thus, equipped with the Weyl product, $W$ is
a deformation quantization of
the formal completion of $\bbC^{2d}$ at the origin. The induced symplectic structure
is the standard one. Let $\frak g$ denote the Lie algebra of continuous derivations of
$W$. Let $\widetilde{\frak g}$ denote the Lie algebra $\displaystyle\frac1{t}W$ with
the bracket $[\displaystyle\frac1{t}f,\displaystyle\frac1{t}g] =
\displaystyle\frac1{t}(f\ast g-g\ast f)$.

The exact sequence
\[
0 @>>> \frac1{t}\bbC[[t]] @>>> \widetilde{\frak g} @>{\frac1{t}ad}>> {\frak g} @>>> 0
\]
is a central extension of Lie algebras.

Let $\deg\hat{z}_i =\deg\hat{z}_i = 1$, $\deg t = 2$, let
$\widetilde{\frak g}_k$ denote the subspace of homogeneous elements of degree $k$,
and let ${\frak g}_k$ denote the image of $\widetilde{\frak g}_k$. Then
$\sum_k\widetilde{\frak g}_k$ (respectively $\sum_k{\frak g}_k$) is a $\bbZ$-graded
dense subalgebra of $\widetilde{\frak g}$ (respectively $\frak g$). 

The Lie algebra ${\frak sp}(2d,\bbC)$ is naturally embedded in both
$\widetilde{\frak g}_0$ and ${\frak g}_0$ as the span of the monomials
$\displaystyle\frac1{t}\hat{z}_i\hat{z}_j$,
$\displaystyle\frac1{t}\hat{z}_i\hat{\xi}_j$, and
$\displaystyle\frac1{t}\hat{\xi}_i\hat{\xi}_j$. The (adjoint) action of
${\frak sp}(2d,\bbC)$ on $\widetilde{\frak g}$ and ${\frak g}$ integrates to an action
of the group $Sp(2d,\bbC)$.

Suppose that $M$ is a complex manifold of dimension $\dim_{\bbC}M = 2d$,
$\omega\in\Gamma(M;\Omega^2_M)$ a symplectic
form on $M$. Let $\cU = \{U_\alpha\}_{\alpha\in I}$ be a cover of $M$ which trivializes
$TM$. Set $U_{\alpha,\beta} = U_\alpha\cap U_\beta$ and let $g_{\alpha,\beta}:
U_{\alpha,\beta} @>>> Sp(2d,\bbC)$ denote the transition functions for $TM$.

Let $\bbW_M$ denote the sheaf of $\cO_M$ algebras which corresponds to the
representation $W$ of $Sp(2d,\bbC)$ (``the sheaf of holomorphic section of the
associated bundle''). A Fedosov connection is a flat $\frak g$-connection $\nabla^\Phi$
on $\bbW_M\otimes_{\cO_M}\cC^\infty_M$ of a particular kind which we proceed to
describe.

Locally (i.e. on a $U_\alpha\in\cU$) a Fedosov connection is given by
$\nabla^\Phi = d + A_\alpha$, where $A_\alpha = \sum_{k\geq -1}A^{(k)}_\alpha$
with $A^{(k)}_\alpha\in
\cA^1_M(U_\alpha)\otimes{\frak g}_k$. Moreover, in Darboux coordinates
$z_1,\dots,z_d,\xi_1\dots\xi_d$, $A^{(-1)}_\alpha =
\displaystyle\frac1{t}\sum_i dz_i\otimes\hat{\xi}_i - d\xi_i\otimes\hat{z}_i$, and
$A^{(0)}_\alpha\in\cA^{1,0}_M(U_\alpha)\otimes{\frak sp}(2d,\bbC)$ so that
$\nabla^{(0)} = d+A^{(0)}_\alpha$ is a torsion free connection.

A Fedosov connection $\nabla^\Phi$ admits a lifting to a
$\widetilde{\frak g}$-connection $\widetilde\nabla^\Phi$ (which is not flat) with
curvature $c(\widetilde\nabla^\Phi)\in\Gamma(M;\cA^{2,cl}_M\widehat\otimes
\displaystyle\frac1{t}\bbC[[t]])$ where $\cA^{2,cl}_M$ is the sheaf of closed 2-forms.

{\em The subsheaf of horizontal sections
$\left(\bbW_M\otimes_{\cO_M}\cC^\infty_M\right)^{\nabla^\Phi}$ is a deformation
quantization of $(M,\omega)$ with characteristic class $\theta$ the cohomology
class of $c(\widetilde\nabla^\Phi)$ which is independent of the choice of a lifting.}

It is shown in \cite{NT3} that every deformation quantization arises in from a
Fedosov connection and that the characteristic classes of isomorphic deformations
coincide.
\subsection{Homological properties of symplectic deformations}
The sheaf of algebras $\bbA^t_M\widehat\boxtimes_{\bbC[[t]]}
\left(\bbA^t_M\right)^{op}$ is contained in a unique symplectic deformation
quantization $\bbA^t_{M^{\times 2}}$ of the symplectic manifold
$(M,\omega)\times(M,-\omega)$.

Let
\[
\left(\bbA^t_M\right)^{\frak{e}}
\overset{def}{=}\delta^{-1}\bbA^t_{M^{\times 2}}\ .
\]
Then there is a
natural map of sheaves of algebras $\bbA^t_M\widehat\otimes_{\bbC[[t]]}
\left(\bbA^t_M\right)^{op} @>>> \left(\bbA^t_M\right)^{\frak{e}}$.

Let
\[
\widehat{C}_p(\bbA^t_M)
\overset{def}{=}
\underset{n}{\varprojlim}
\left(\bbA^t_M\right)^{\widehat\otimes_{\bbC[[t]]}p+1}/\cI_{(p+1)}^n\ ,
\]
where (the two-sided ideal) $\cI_{(p+1)}$ is the kernel of the composition
of the multiplication map with the symbol map
$\left(\bbA^t_M\right)^{\widehat\otimes_{\bbC[[t]]}p+1}
@>>> \bbA^t_M @>>> \cO_M$. The Hochschild differential \eqref{diffl:b}
extends to this setting and gives rise to the complex $\widehat{C}_\bullet
(\bbA^t_M)$ and the natural map of complexes
$C_\bullet(\bbA^t_M) @>>> \widehat{C}_\bullet(\bbA^t_M)$.
The complex $\widehat{C}_\bullet(\bbA^t_M)$ represents
$\bbA^t_M\otimes^{\bbL}_{\left(\bbA^t_M\right)^{\frak{e}}}
\bbA^t_M$ in the derived category.

The cyclic complexes $\widehat{CC}^{-}_\bullet(\bbA^t_M)$ and
$\widehat{CC}^{per}_\bullet(\bbA^t_M)$ are defined as in
\ref{subsection:C-CC} starting from $\widehat{C}_p(\bbA^t_M)$
The analogs of \eqref{inclusions} and \eqref{ses:CC-C} hold and there are
natural maps of complexes
$CC^{-}_\bullet(\bbA^t_M) @>>> 
\widehat{CC}^{-}_\bullet(\bbA^t_M)$ and
$CC^{per}_\bullet(\bbA^t_M) @>>> 
\widehat{CC}^{per}_\bullet(\bbA^t_M)$.

\begin{prop}\label{prop:CCA}
Suppose that $\bbA^t_M$ is a symplectic deformation quantization of $M$.
\begin{enumerate}
\item Let $x_1,\ldots,x_d,\xi_1,\ldots,\xi_d$ be a local Darboux coordinate
system on $M$. Then the local section $\Phi^{\bbA}$ of
$H^{-2d}\widehat{C}_\bullet(\bbA^t_M)[t^{-1}]$
represented by $Alt(1\otimes x_1\otimes\cdots\otimes x_d\otimes
\displaystyle\frac{\xi_1}{t}\otimes\cdots\otimes\displaystyle\frac{\xi_d}
{t})$
is independent of the choice of the Darboux coordinates.

\item The section $\Phi^{\bbA}$ generates
$H^{-2d}\widehat{C}_\bullet(\bbA^t_M)[t^{-1}]$ as a $\bbC_M[t^{-1},t]]$-module
and $H^{p}\widehat{C}_\bullet(\bbA^t_M)[t^{-1}]=0$ for $p\neq -2d$, thus
and there is a unique isomorphism
\begin{equation}\label{map:muA}
\mu_{\bbA}:\widehat{C}_\bullet(\bbA^t_M)[t^{-1}]
@>{\isomo}>>
\bbC_M[t^{-1},t]][2d]
\end{equation}
in the derived category of sheaves which maps $\Phi^{\bbA}$ to $1$.

\item There exists an isomorphism of filtered complexes
\[
(\widehat{CC}^{-}_\bullet(\bbA^t_M)[t^{-1}],b+uB,F_\bullet)
\isomo
(\widehat{C}_\bullet(\bbA^t_M)[t^{-1}][[u]],b,F_\bullet)
\]
which induces the identity map on the associated graded objects and any two such
are homotopic.

\item There exists an isomorphism of filtered complexes
\[
(\widehat{CC}^{per}_\bullet(\bbA^t_M)[t^{-1}],b+uB,F_\bullet)
\isomo
(\widehat{C}_\bullet(\bbA^t_M)[t^{-1}][u^{-1},u]],b,F_\bullet)
\]
which induces the identity map on the associated graded objects and any two such
are homotopic.
\end{enumerate}
\end{prop}
\begin{pf}
Since the statement is local it is sufficient to consider the stalk at the
origin of the standard deformation of $\bbC^{2d}$. Note that the completion
of $\bbA^t_{{\bbC^{2d}},0}$ with respect to the powers of the two-sided
ideal which is the kernel of the composition $\bbA^t_{{\bbC^{2d}},0}
@>{\sigma}>> \cO_{{\bbC^{2d}},0} @>>> \bbC_0$ coincides with the Weyl algebra
$W=W(\bbC^{2d})$. The induced map $\widehat{C}_\bullet(\bbA^t_{\bbC^{2d}})_0
@>>> \widehat{C}_\bullet(W)$ is a quasiisomorphism and the calculation
of Hochschild homology follows from the standard results for the Weyl algebra.
To prove the last two statements one argues as in Lemma \ref{lemma:CC-E} and
Corollary \ref{cor:CC-E}.
\end{pf}

By abuse of notation we will denote by $\mu_{\bbA}$ the composition
\[
C_\bullet(\bbA^t_M) @>>> \widehat{C}_\bullet(\bbA^t_M)[t^{-1}]
@>{\text{\eqref{map:muA}}}>>
\bbC_M[t^{-1},t]][2d]
\]
as well as \eqref{map:muA}.

Proceeding as in \ref{subsec:trace-density} and using Proposition \ref{prop:CCA}
one defines the morphism
\[
\tilde\mu_{\bbA} : CC^{per}_\bullet(\bbA^t_M)[t^{-1}]
@>>> \bbC_M[t^{-1},t]][2d][u^{-1},u]]\ .
\]

\subsection{Deformation quantization of a cotangent bundle}
For applications of our results to microlocal analysis we will require the
particular deformation quantization of a cotangent bundle furnished by the
following proposition.

\begin{prop}\label{prop:DQcotan}
There exists a deformation quantization $\bbA^t_{T^*X}$ of $T^*X$
and faithfully flat map $\cR\cE_X @>>> \bbA^t_{T^*X}$
of algebras over $\bbC [t ]$. The characteristic class $\theta$ of the 
deformation $\bbA^t_{T^*X}$ is equal to
$\displaystyle\frac12\pi^*c_1(TX)$ (note that the symplectic form is exact in
this case).
\end{prop}
\begin{pf}
Suppose that $x_1,\ldots,x_d$ is a local coordinate system on $X$ and let
$\xi_i$ denote the symbol of $\displaystyle\frac{\partial}{\partial x_i}$
viewed as a function on $T^*X$. Then $x_1,\ldots,x_d,\xi_1,\ldots,\xi_d$
form a Darboux coordinate system on $T^*X$. The map alluded to in Proposition
\ref{prop:DQcotan} is given in local coordinates by $x_i\mapsto x_i$,
$t\displaystyle\frac{\partial}{\partial x_i}\mapsto\xi_i$.

The only non-trivial part is the calculation of the characteristic class $\theta$.
This will be achieved by explicitly manufacturing a Fedosov connection
$\nabla^\Phi$ on the sheaf of Weyl algebras $\bbW_{T^*X}$ whose sheaf of horizontal
sections is $\bbA^t_{T^*X}$ and computing the Weyl curvature.

Let $\cU = \{ U_\alpha \}_{\alpha\in I}$ denote a cover of X which trivializes $TX$.
Set $U_{\alpha,\beta} = U_\alpha\cap U_\beta$ and let
$g_{\alpha,\beta} : U_{\alpha,\beta} @>>> GL(d,\bbC)$ denote the transition functions
for $TX$.

Let ${\frak w}_d$ denote the Lie algebra of continuous derivations of the (topological)
$\bbC$-algebra $\bbC[[\hat{z}_1,\dots,\hat{z}_d]]$. Let ${\frak w}_{d,k} =
\{\sum_i P_i\displaystyle\frac{\partial\ }{\partial\hat{z}_i}\}$ where $P_i$ are
homogeneous polynomials (in $\hat{z}_1,\dots,\hat{z}_d$) of degree $\deg P_i = k+1$.
Then $\left[{\frak w}_{d,k},{\frak w}_{d,l}\right]\subseteq {\frak w}_{d,k+l}$ and
$\sum_k{\frak w}_{d,k}$ is a dence subalgebra of ${\frak w}_d$. The Lie algebra
${\frak gl}(d,\bbC)$ is identified with ${\frak w}_{d,0}$ by the map
$\left(a_{i,j}\right)\mapsto\sum_{i,j}a_{i,j}\hat{z}_i\displaystyle
\frac{\partial\ }{\partial\hat{z}_j}$. The (adjoint) action of ${\frak gl}(d,\bbC)$
on ${\frak w}_d$ integrates to the action of $GL(d,\bbC)$. The action of
${\frak w}_{d,0}$ on ${\frak w}_{d,-1}$ is identified with the standard representation
of ${\frak gl}(d,\bbC)$ on $\bbC^d$.

We begin with a torsion free ${\frak gl}(d,\bbC)$-connection $\nabla^{(0)}$ on $TX$
given on $U_\alpha$ by $d + A^{(0)}_\alpha$, where $A^{(0)}_\alpha\in\cA^{1,0}_X\otimes
{\frak gl}(d,\bbC)$. The connection $\nabla^{(0)}$ extends to a flat
${\frak w}_d$-connection (called a Kazhdan connection) $\nabla^K$ which is, locally,
of the form $d + A_\alpha$, where $A_\alpha\in\cA^1_X\widehat\otimes{\frak w}_d$,
$A_\alpha = \sum_k A^{(k)}$ with
$A^{(-1)}_\alpha=-\sum_i dz_i\otimes\displaystyle\frac{\partial\ }{\partial\hat{z}_i}$
and $A^{(0)}_\alpha$ as above.

Let $\widehat S = \bbC[[\hat{z}_1,\dots,\hat{z}_d]]$; this is a
$({\frak w}_d,GL(d,\bbC))$-module. Let $\widehat\cS$ denote the corresponding
sheaf of $\cO_X$-algebras (``holomorphic sections of the associated bundle''). The
Kazhdan connection $\nabla^K$ is a connection on
$\widehat\cS\otimes_{\cO_X}\cC^\infty_X$, and the sheaf of
horizontal sections $\left(\widehat\cS\otimes_{\cO_X}\cC^\infty_X\right)^{\nabla^K}$
is a sheaf of algebras (since ${\frak w}_d$ acts by derivations) isomorphic to $\cO_X$.

Let $\widehat R = \bbC[[\hat{z}_1,\dots,\hat{z}_d,\hat{\xi}_1,\dots,\hat{\xi}_d,t]]$.
We consider $\widehat R$ as an algebra with the Weyl product denoted $\ast$. The
algebra $\widehat R$ is an algebra over $\widehat S$ in the obvious way and is
isomorphic to the completion of the Rees algebra of the algebra
of differential operators on $\widehat S$ by the assignment
$t\displaystyle\frac{\partial\ }{\partial\hat{z}_i}\mapsto\hat{\xi}_i$.

Let $Der(\widehat R)$ denote the Lie algebra of continuous, $\bbC[[t]]$-linear 
derivations of $\widehat R$. The exact sequence
\[
0 @>>> \frac1{t}\bbC[[t]] @>>> \frac1{t}\widehat R
@>{\frac1{t}ad}>> Der(\widehat R) @>>> 0
\]
is a central extension of Lie algebras. Just as in the case of the Weyl algebra
in \ref{subsection:Fedosoff},
the Lie algebra ${\frak sp}(2d,\bbC)$ is contained in $\frac1{t}\widehat R$ and
$ Der(\widehat R)$ is the standard way. The inclusion of ${\frak gl}(d,\bbC)$ into
${\frak sp}(2d,\bbC)$ is given in terms of the standard embedding by
$(a_{i,j})\mapsto\sum_{i,j}a_{i,j}\hat{z}_i\hat{\xi}_j$. The action of
${\frak gl}(d,\bbC)$ integrates to an action of $GL(d,\bbC)$.

The identification with the differential
operators (above) restricts to the map of Lie algebras
$i : {\frak w}_d @>>>\displaystyle\frac1{t}\widehat R$ whose restriction to
${\frak gl}(d,\bbC)$ is
\begin{equation}\label{map:gl-to-R}
\left(a_{i,j}\right)\mapsto
\sum_{i,j}a_{i,j}\hat{z}_i\ast\displaystyle\frac{\hat{\xi}_j}{t}
= \sum_{i,j}a_{i,j}\hat{z}_i\displaystyle\frac{\hat{\xi}_j}{t} - 
\frac12 tr\left(a_{i,j}\right)\ .
\end{equation}
Note that the first summand in the last expression is the standard embedding
of ${\frak gl}(d,\bbC)$. Since the second summand in the last expression is
central the standard embedding and $i$ coincide in $Der(\widehat R)$.

Let $\widehat\cR$ denote the sheaf of $\cO_X$-algebras on $X$ which corresponds to
$\widehat R$. 
The Kazhdan connection $\nabla^K$ gives rise to a connection, still denoted
$\nabla^K$, on $\widehat\cR\otimes_{\cO_X}\cC^\infty_X$ via the composition
$\displaystyle\frac1{t}ad\circ i$. The sheaf of horizontal sections
$\left(\widehat\cR\otimes_{\cO_X}\cC^\infty_X\right)^{\nabla^K}$ is isomorphic to
the completion of the Rees algebra $\cR\cD_X$.

There is a canonical
isomorphism of algeras $\widehat\cR\isomo\bbW_{T^*X}\otimes_{\cO_{T^*X}}\cO_X$.
Since the sheaves of algebras $\bbW_{T^*X}\otimes_{\cO_{T^*X}}\cC^\infty_{T^*X}$ and
$\pi^*\widehat\cR\otimes_{\cO_{T^*X}}\cC^\infty_{T^*X}$ are (non-canonically) isomorphic,
it suffices to construct a ``Fedosov'' $Der(\widehat R)$-connection on the latter and
this is what we will do explicitly, at least in the formal neighborhood of the zero
section of $T^*X$.

Let $z_1,\dots,z_d$ denote local coordinates on $X$, $\xi_i$ the symbol of
$\displaystyle\frac{\partial\ }{\partial z_i}$, so that $z_1,\dots,z_d,\xi_1\dots\xi_d$
is a local Darboux coordinate system on $T^*X$. The automorphism $\Psi$ of
sheaf of algebras $\pi^*\cR$ given in local coordinates by the formula
$\Psi=\exp(ad(-\sum_i\xi_i\displaystyle\frac{\hat{z}_i}{t}))$ is, in fact globally well
defined on the formal neighborhood of the zero section of $T^*X$. The $\Psi$-conjugate
of the pull-back $\pi^*\nabla^K$ of the Kazhdan connection is the desired Fedosov
connection: $\nabla^\Phi = \left(\nabla^K\right)^\Psi$.

Suppose $\tilde{\pi^*\nabla^K}$ is a lifting of $\pi^*\nabla^K$ to a
$\displaystyle\frac1{t}\widehat R$-connection. Then
$\left(\tilde{\pi^*\nabla^K}\right)^\Psi$ is a lifting of $\nabla^\Phi$, and
$c(\tilde{\pi^*\nabla^K}) = 
c\left(\left(\tilde{\pi^*\nabla^K}\right)^\Psi\right)$ since
the curvature is central and $\Psi$ is an ``inner'' automorphism.

In order to find an explicit lifting $\tilde{\pi^*\nabla^K}$, consider the
$\displaystyle\frac1{t}\widehat R$-valued forms $i(A_\alpha)$. Since the
(${\frak w}_d$-valued) connection forms $A_\alpha$ satisfy
\[
A_\alpha = dg^t_{\alpha,\beta}(g^t_{\alpha,\beta})^{-1}
+ g^t_{\alpha,\beta}A_\beta(g^t_{\alpha,\beta})^{-1}\ ,
\]
the forms $iA_\alpha$ satisfy
\begin{eqnarray*}
i(A_\alpha) & = & i(dg^t_{\alpha,\beta}(g^t_{\alpha,\beta})^{-1}) +
		  i(g^t_{\alpha,\beta}A_\beta(g^t_{\alpha,\beta})^{-1}) \\
& & = dg^t_{\alpha,\beta}(g^t_{\alpha,\beta})^{-1} -
		\frac12 tr(dg^t_{\alpha,\beta}(g^t_{\alpha,\beta})^{-1}) +
		g^t_{\alpha,\beta}i(A_\beta)(g^t_{\alpha,\beta})^{-1}\ ,
\end{eqnarray*}
the last equality obtained using \eqref{map:gl-to-R}.
Since $A^{(0)}_\alpha$ are (${\frak gl}(d,\bbC)$-valued) connection forms and trace
remains invariant under conjugation the identity
\[
\frac12 tr(A^{(0)}_\alpha) = \frac12 tr(dg^t_{\alpha,\beta}(g^t_{\alpha,\beta})^{-1})
+ \frac12 tr(A^{(0)}_\beta)
\]
holds. Adding the latter identity to the preceding one one easily concludes that
the collection of $\displaystyle\frac1{t}\widehat R$-valued forms $\widetilde A_\alpha
= i(A_\alpha) + \displaystyle\frac12 tr(A^{(0)}_\alpha)$ determine a connection
$\tilde\nabla^K$ which lifts the Kazhdan connection.

Since the Kazhdan connection is flat and $i$ is a morphism of Lie algebras, it follows
that $c(\tilde\nabla) = \displaystyle\frac12 tr((\nabla^{(0)})^2)$ which represents
$\displaystyle\frac12 c_1(TX)$ in de Rham cohomology.

Hence, the cohomology class of the curvature of the connection $\tilde\pi^*\nabla^K$
is equal to $\displaystyle\frac12 \pi^*c_1(TX)$ and so is the characteristic class of
the deformation.
\end{pf}

In what follows $\bbA^t_{T^*X}$ will always refer to the deformation furnished by
Proposition \ref{prop:DQcotan}.

\begin{lemma}
The diagram
\[
\begin{CD}
C_\bullet(\cR\cE_X)[t^{-1}] @>{\isomo}>> C_\bullet(\cE_X)[t^{-1},t]\\
@VVV				@VV{\mu_{\cE}}V \\
C_\bullet(\bbA^t_{T^*X})[t^{-1}]
@>{\mu^t_{\bbA}}>> \bbC_{T^*X}[t^{-1},t]][2d]
\end{CD}
\]
is commutative.
\end{lemma}
\begin{pf}
The statement is equivalent to the commutativity of the diagram of respective
cohomology sheaves in degree $-2d$. The latter fact follows from the local
coordinate representation
of $\Phi^{\cE}$, $\Phi^{\bbA}$, the map $\cR\cE_X @>>> \bbA^t_{T^*X}$,
and the definitions of $\mu_{\cE}$ and $\mu^t_{\bbA}$.
\end{pf}
\begin{cor}\label{cor:EtoA}
The diagram
\[
\begin{CD}
CC_\bullet(\cR\cE_X)[t^{-1}]
@>{\isomo}>> 
CC_\bullet(\cE_X)[t^{-1},t]
\\
@VVV				@VV{\tilde\mu_{\cE}}V 
\\
CC_\bullet(\bbA^t_{T^*X})[t^{-1}]
@>{\tilde\mu^t_{\bbA}}>>
\bbC_{T^*X}[t^{-1},t]][2d][u^{-1},u]]
\end{CD}
\]
is commutative.
\end{cor}

\subsection{The Riemann-Roch theorem for periodic cyclic cocycles}
We are finally ready to state the Riemann-Roch type theorem for symplectic
deformations. Consider a complex manifold $M$ and a symplectic deformation
quantization $\bbA^t_M$ of $M$ with the characteristic class $\theta$.

The deformation quantization $\bbA^t_M$ induces the symplectic
structure on $M$. In particular the symplectic structure determines
a reduction of the structure group of $TM$ to $Sp(2d,\bbC)$. The characteristic
classes of $Sp(2d,\bbC)$-bundles are given by the invariant functions on the
Lie algebra ${\frak sp}(2d,\bbC)$. In particular the class $\widehat A(TM)$
corresponds to the function
\[
{\frak sp}(2d,\bbC)\ni X\mapsto \det\left(\frac{\ad(X)}
{\exp(\ad(\frac{X}2)) - \exp(\ad(-\frac{X}2))}\right)\ .
\]

\begin{thm}\label{thm:main}
The diagram
\[
\begin{CD}
CC^{per}_\bullet(\bbA^t_M) @>{\sigma}>> CC^{per}_\bullet(\cO_M) \\
@V{\iota}VV @VV
{\tilde\mu_{\cO}\smile \widehat A(TM)\smile e^\theta}V \\
CC^{per}_\bullet(\bbA^t_M)[t^{-1}] @>{\tilde\mu^t_{\bbA}}>>
\bbC_M[t^{-1},t]][u^{-1},u]]
\end{CD}
\]
is commutative.
\end{thm}

Theorem \ref{thm:main} may be deduced from \cite{NT1} and \cite{NT2}.
A complete proof will appear in the sequel to the present paper.

\section{Appendix: Review of Hochschild and cyclic homology}
\label{section:CC}
In this section we review the basic definitions of the Hochschild and the
(negative, periodic) cyclic complex of an algebra, of a sheaf of algebras
on a space as well as of the ``topological'' versions of the above. In
particular we establish notational conventions with regard to Hochschild and
cyclic complexes which are used in the body of the paper.

In addition we review the results concerning the Hochschild and cyclic homologies
of certain examples of (sheaves of) algebras which appear in this paper for
the reader's convenience.

\subsection{Hohschild and cyclic complexes of algebras}
\label{subsection:C-CC}
Let $k$ denote a commutative algebra over a field of characteristic zero
and let $A$ be a flat $k$-algebra with $1_A\cdot k$ contained in the center,
not necessarily commutative. Let $\overline A = A/k$, and let
$C_p(A)\overset{def}{=}A\otimes_k\overline A^{\otimes_k p}$ and let
\begin{eqnarray}\label{diffl:b}
b : C_p(A) & @>>> & C_{p-1}(A) \\ \nonumber
a_0\otimes\cdots\otimes a_p & \mapsto & 
(-1)^p a_pa_0\otimes\cdots\otimes a_{p-1} + \\ & &
\sum_{i=0}^{p-1}(-1)^ia_0\otimes\cdots\otimes a_ia_{i+1}\otimes\cdots\otimes
a_p\ .\nonumber
\end{eqnarray}
Then $b^2 = 0$ and the complex $(C_\bullet(A), b)$, called {\em the standard
Hochschild complex of $A$} represents $A\otimes^{\bbL}_{A\otimes_kA^{op}}A$
in the derived category of $k$-modules.

The map
\begin{eqnarray}\label{diffl:B}
B : C_p(A) & @>>> & C_{p+1}(A) \\ \nonumber
a_0\otimes\cdots\otimes a_p & \mapsto & \sum_{i=0}^p (-1)^{pi}
1\otimes a_i\otimes\cdots\otimes a_p\otimes a_0\otimes\cdots\otimes a_{i-1}
\end{eqnarray}
satisfies $B^2 = 0$ and $[B,b] = 0$ and therefore defines a map of complexes
\[
B: C_\bullet(A) @>>> C_\bullet(A)[-1]\ .
\]
For $i,j,p\in\bbZ$ let
\begin{eqnarray*}
CC^-_p(A) & = & \prod_{\overset{i\geq 0}{i+j = p \mod 2}} C_{i+j}(A) \\
CC^{per}_p(A) & = & \prod_{i+j = p \mod 2} C_{i+j}(A)\ .
\end{eqnarray*}
The complex $(CC^-_\bullet(A),B+b)$ (respectively
$(CC^{per}_\bullet(A),B+b)$) is called the {\em negative} (respectively
{\em periodic}) {\em cyclic complex of $A$}.

There are inclusions of complexes
\begin{equation}\label{inclusions}
CC^-_\bullet(A)[-2]\hookrightarrow CC^-_\bullet(A)\hookrightarrow
CC^{per}_\bullet(A)
\end{equation}
and the short exact sequence
\begin{equation}\label{ses:CC-C}
0 @>>> CC^-_\bullet(A)[-2] @>>> CC^-_\bullet(A) @>>> C_\bullet(A) @>>> 0\ .
\end{equation}

In what follows we will use the notation of Getzler and Jones (\cite{GJ}).
Let $u$ denote a variable of degree $-2$ (with respect to the homological
grading. Then the negative and periodic
cyclic complexes are described by the following formulas:
\begin{eqnarray} \label{ses:gejo}
CC^-_\bullet(A) & = & (C_{\bullet}(A)[[u]], b + uB) \\
\label{ses:gejoper}
CC^{per}_\bullet(A) & = & (C_{\bullet}(A)[[u, u^{-1}], b + uB)\ .
\end{eqnarray}

\subsection{Hochschild and cyclic complexes of sheaves of algebras}
Suppose that $X$ is a topological space and $\cA$ is a flat sheaf of 
$k$-algebras on $X$ such that there is a global section $1\in\Gamma(X;\cA)$ 
which
restricts to $1_{\cA_x}$ and $1_{\cA_x}\cdot k$ is contained in the center of 
$\cA_x$ for every point $x\in X$. Let $C_\bullet(\cA)$ (respectively
$CC^-_\bullet(\cA),\ CC^{per}_\bullet(\cA)$) denote the complex
of sheaves of $k$-modules associated to the presheaf with value
$C_\bullet(\cA(U))$ (respectively
$CC^-_\bullet(\cA(U)),\ CC^{per}_\bullet(\cA(U))$) on an open subset $U$ of
$X$. Then $C_\bullet(\cA)$ represents
$\cA\otimes^{\bbL}_{\cA\otimes_k\cA^{op}}\cA$ in the derived category of
sheaves of $k$-modules on $X$.

\end{document}